\documentclass[12pt]{amsart}
\usepackage{amscd}
\newtheorem{thm}{Theorem}
\newtheorem{prop}{Proposition}
\newtheorem{lem}{Lemma}

\theoremstyle{remark}
\newtheorem{rem}{Remark}

\theoremstyle{definition}
\newtheorem{defn}{Definition}

\newcommand{\lra}{\longrightarrow}

\newcommand{\R}{\mathbb{ R}}

\newcommand{\rk}{\operatorname{rk}}

\DeclareMathOperator{\dummygg}{\mathfrak{g}}
\renewcommand{\gg}{\dummygg}
\DeclareMathOperator{\hh}{\mathfrak{h}}
\DeclareMathOperator{\TT}{\mathfrak{t}}
\DeclareMathOperator{\s}{\mathfrak{s}}

\title{On formality of generalised symmetric spaces}
\author{D.~Kotschick}
\address{Mathematisches Institut, Ludwig-Maximilians-Universit\"at
M\"unchen, Theresienstr.~39, 80333 M\"unchen, Germany}
\email{dieter@member.ams.org}
\author{S.~Terzi\'c}
\address{Mathematisches Institut, Ludwig-Maximilians-Universit\"at
M\"unchen, Theresienstr.~39, 80333 M\"unchen, Germany}
\email{terzic@rz.mathematik.uni-muenchen.de}
\email{sterzic@cg.ac.yu}
\thanks{The authors are members of the {\sl European Differential Geometry 
Endeavour} (EDGE), Research Training Network HPRN-CT-2000-00101, 
supported by The European Human Potential Programme; the second author is 
supported by the {\it DFG Graduiertenkolleg ``Mathematik im Bereich ihrer
Wechselwirkung mit der Physik''}.}
\date{\today; MSC2000: primary 57T15, secondary 53C35, 55P62, 58A14}

\begin{document}

\begin{abstract}
We prove that all generalised symmetric spaces of compact simple Lie 
groups are formal in the sense of Sullivan. Nevertheless, many of 
them, including all the non-symmetric flag manifolds, do not admit 
Riemannian metrics for which all products of harmonic forms are 
harmonic.
\end{abstract}

\maketitle

\section{Introduction}\label{s:intro}

In this paper we discuss formality properties of certain compact
homogeneous spaces $G/H$, with $G$ a compact connected Lie group and
$H$ a closed subgroup. 
We shall discuss formality in the sense of Sullivan's rational 
homotopy theory~\cite{tokyo} and geometric formality in the sense 
of~\cite{K}. 

There are some classes of compact homogeneous spaces which are
well-known to be formal in the sense of Sullivan, for example
the symmetric spaces and those homogeneous spaces with $\rk G=\rk H$.
We shall see that it is an immediate consequence of earlier work of 
the second author~\cite{T} that in fact all \emph{generalised}
symmetric spaces\footnote{These spaces are sometimes called 
\emph{$k$--symmetric}.} of compact simple Lie groups are formal in the 
sense of Sullivan.

The notion of geometric formality was introduced by the first author
in~\cite{K}. A smooth manifold is said to be geometrically formal if
it admits a Riemannian metric for which all wedge products of 
harmonic forms are harmonic. This clearly implies formality in the 
sense of Sullivan, and is even more restrictive. As compact
symmetric spaces are the classical examples of geometrically formal 
manifolds, it is natural to explore this notion in the context of 
generalised symmetric spaces. Although these turn out to be formal 
in the sense of Sullivan and also satisfy all the restrictions on 
geometrically formal manifolds found in~\cite{K}, we shall prove 
here that many of them are \emph{not} geometrically formal. 
Some of our examples have $\rk G = \rk H$, whereas others do not. 

At the time of writing it remains unclear whether there is a
reasonable class of non-symmetric compact homogeneous spaces which
are geometrically formal.

In Section~\ref{s:gsscoh} we collect some classical results on the 
cohomology of compact homogeneous spaces, and we summarise the results
we shall need from~\cite{T0,T} on the classification of generalised 
symmetric spaces and their cohomology. These results are used in 
Section~\ref{s:formal} to conclude that all generalised
symmetric spaces of compact simple Lie groups are formal in the 
sense of Sullivan. Section~\ref{s:flag} makes explicit the additive 
generators and multiplicative relations between them for the 
cohomology algebras of the flag manifolds. This is then used in 
Section~\ref{s:geoform} to prove that the non-symmetric flag manifolds and 
several other classes of generalised symmetric spaces are not 
geometrically formal.

\section{The real cohomology of compact homogeneous spaces}\label{s:gsscoh}

Let $G$ be compact connected Lie group, and $H\subset G$ a connected closed
subgroup. We denote by $\TT$ and $\s$ the maximal abelian subalgebras
of the Lie algebras $\gg$ and $\hh$ respectively, and by $BG$ the
classifying space of $G$.

By the Hopf theorem~\cite{Borel}, $H^{*}(G)$ is an exterior algebra
on universal transgressive elements $z_{1},\ldots, z_{n}$. The 
Cartan-Chevalley theorem~\cite{Borel} implies that $H^{*}(BG)$ is the 
ring of $W$-invariant polynomials on $\TT$ with real coefficients,
where $W$ is the Weyl group of $\gg$ relative to $\TT$. For all compact
simple Lie groups the generators of the Weyl invariants are 
well-known~\cite{Doan Kuin'}.
Coordinates $x_1,\ldots ,x_n$ on $\TT$ expressing the Weyl invariant polynomials in 
classical form will be called canonical coordinates.
Let $y_{1},\ldots, y_{n}$ correspond to $z_{1}, \ldots, z_{n}$ by
transgression in the universal $G$-bundle over $BG$. Then $H^{*}(BG)$ 
is generated by $y_{1}, \ldots, y_{n}$~\cite{Borel}.

We consider the map $\rho^{*}(H, G)\colon{\R}[\TT]^{W_G}\rightarrow
{\R}[\s]^{W_H}$ assigning to each polynomial in ${\R}[\TT]^{W_G}$ its
restriction to $\s$.
The {\em Cartan algebra} of the homogeneous space $G/H$ is the algebra
$C = {\R}[\s]^{W_H}\otimes H^{*}(G)$ endowed with the following differential
$d$:
\[ 
d(1\otimes z_{i}) = \rho ^{*}(H, G)y_{i}\otimes {1} \; \; \; \; \; (1\leq i\leq n),
\]
\[ 
d(b\otimes 1) = 0 \; \;\mbox{for} \; \; b\in {\R}[\s]^{W_H} \ .
\]
The name derives from the following celebrated result:
\begin{thm}[Cartan]
The real cohomology algebra of the homogeneous space $G/H$ is isomorphic to
the cohomology algebra of its Cartan algebra $(C,d)$.
\end{thm}
This theorem in principle computes the cohomology of $G/H$. In
practice, however, one still needs information about the map $\rho
^{*}(H,G)$
in order to obtain an explicit result.

Before giving a summary of the calculations for generalised symmetric 
space that we shall need, we recall some earlier applications of 
Cartan's theorem.

\subsection{Homogeneous spaces with $\rk G = \rk H$}
In his classical paper~\cite{Borel}, Borel studied the cohomology rings 
of homogeneous spaces with $\rk G = \rk H$. For these he showed that 
$\rho ^{*}(H,G)$ is injective, which implies:
\begin{thm}[\cite{Borel}]\label{first}
The real cohomology algebra of a compact homogeneous space $G/H$ with
$\rk G = \rk H$ is given by
\[ 
H^{*}(G/H)\cong {\R}[\TT]^{W_H}/\langle {\R}[\TT]^{W_G}\rangle^{+} \ ,
\]
where $\langle {\R}[\TT]^{W_G}\rangle^{+}$ is the ideal in 
${\R}[\TT]^{W_H}$ generated by the elements of ${\R}[\TT]^{W_G}$ of 
positive degrees.
\end{thm}
Note that, in order to obtain a more explicit formula, one also
requires information about the transition functions
between canonical coordinates of the group $G$ and its subgroup $H$.

\subsection{Symmetric spaces}\label{ss:symm}
Borel~\cite{Borel} also calculated the cohomology algebras of the symmetric
spaces 
$SU(n)/SO(n)$ and $SU(2m)/Sp(m)$, which have $\rk H < \rk G$. The cohomology
algebras of the remaining three kinds of symmetric spaces,
$SO(2l)/(SO(2m+1)\times SO(2l-2m-1))$, $E_6/F_4$ and $E_6/PSp(4)$,
were calculated by Takeuchi~\cite{Takeuchi}.

\subsection{Homogeneous spaces of Cartan type}
There is a larger class of compact homogeneous spaces for which Cartan's
theorem can be used directly, which we call {\em homogeneous spaces of 
Cartan type}. They are called {\em normal position homogeneous
spaces} in~\cite{Doan Kuin'}, and {\em Cartan 
pairs} $(G,H)$ in~\cite{GHV}. We say that the 
homogeneous space $G/H$ with $\rk G = n$ and $\rk H = k$ is 
\emph{of Cartan type} if one can choose generators
$P_1, \ldots, P_n$ of ${\R}[\TT]^{W_G}$ in such a way that 
$\rho^{*}(H,G)P_{k+1},\ldots,\rho ^{*}(H, G)P_{n}$ belong to the ideal 
in ${\R}[\s]^{W_H}$ generated by 
$\rho ^{*}(H, G)P_{1},\ldots,\rho ^{*}(H, G)P_{k}$.
The following theorem is proved in~\cite{GHV,Onishchik}:
\begin{thm}\label{normal}
Let $G/H$ be a compact homogeneous space of Cartan type  with 
$\rk G = n$ and $\rk H = k$. Then its cohomology algebra is given by
\begin{equation}\label{eq:Cartan}
H^{*}(G/H)\cong {\R}[\s]^{W_H}/\langle\rho^{*}(H,G)({\R}[\TT]^{W_G})\rangle
\otimes\wedge (z_{k+1},\ldots,z_{n}) \ ,
\end{equation}
where the $z_{i}$ are universal transgressive generators of $H^{*}(G)$.
\end{thm}
We shall refer to the first and second factors in~\eqref{eq:Cartan}
as the {\it polynomial} and the {\it exterior algebra} parts of the 
cohomology.

Note that deciding whether a homogeneous space $G/H$ is of Cartan type, 
or not, is almost equivalent to calculating the map 
$\rho^{*} (H,G)$, and is therefore quite difficult in general.

\begin{rem}
The Poincar\'e polynomial for a homogeneous space $G/H$ of Cartan type is 
given by 
\begin{equation}\label{poincare}
p(G/H, t) = \prod _{i=1}^{k}\frac{1-t^{2k_i}}{1-t^{2l_i}}\prod _{i=k+1}^{n}
(1+t^{2k_{i}-1}) \ ,
\end{equation}
where $k_i$ $(1\leq i\leq n)$ are the exponents of $G$ and $l_i$ 
$(1\leq i\leq k)$ are the exponents of the subgroup $H$. Compare~\cite{GHV}.
\end{rem}
The following lemma provides an important fact about
fibrations between homogeneous spaces.
\begin{lem}\label{l:colapses}
Let $G/H$ and $G/L$ be Cartan type homogeneous spaces with the same 
exterior algebra parts of their cohomologies, and such that $H\subset L$. 
Then the restriction to the fiber $L/H$ of the fibration 
$G/H\lra G/L$ is surjective in real cohomology.
\end{lem}
\begin{proof}
Since the spaces $G/L$ and $G/H$ have the same 
exterior algebra parts of their cohomologies, obviously $\rk L=\rk H$
and from~\eqref{poincare} it follows that 
\[
p(G/H, t) = p(G/L, t)\cdot p(L/H, t) \ .
\]
Thus, the spectral sequence of the fibration collapses.
Now the Leray-Hirsch theorem implies that the restriction to the fiber 
is surjective in real cohomology.
\end{proof}

\subsection{Generalised symmetric spaces}\label{s:gss}

There are several ways of generalising the notion of a symmetric space.
The spaces we shall consider here have been studied by many authors, 
see e.~g.~\cite{G,Kow,V-F,W-G}. They are sometimes called $k$--symmetric,
where $k$ is the order, which we prefer to denote by $m$ below.

\begin{defn}
A generalised symmetric space
of order $m$ is a triple $(G, H, \Theta)$, where $G$ is
a connected Lie group, $H\subset G$ is a closed subgroup, and $\Theta$ is
an automorphism of finite order $m$ of the group $G$ satisfying
\[ 
G_{0}^{\Theta}\subseteq H\subseteq G^{\Theta} \ ,
\]
where $G^{\Theta}$ is the fixed point set of $\Theta$ and $G_{0}^{\Theta}$
is its identity component.
\end{defn}
Obviously, for $m = 2$ these are the usual symmetric spaces.
The ``space'' underlying a generalised symmetric space is the
homogeneous space $G/H$.
Just as in the case of symmetric spaces, generalised symmetric spaces
of order $m$ in the sense of the above definition can be characterised as
Riemannian manifolds admitting certain geodesic symmetries of order
$m$, see for example~\cite{G}. 

The class of generalised symmetric spaces is a lot larger than that
of symmetric spaces; it is easy to see that many generalised 
symmetric spaces do not have the homotopy type of any symmetric space.

For semi-simple and simply connected Lie groups $G$ all fixed point 
subgroups are connected, and there is a bijection between generalised
symmetric spaces and generalised symmetric algebras~\cite{Vinberg}.
One can then discuss triples $(\gg, \gg^{\Theta}, \Theta)$, for 
simple Lie algebras $\gg$ of compact Lie groups with a finite order
automorphism $\Theta$. Assuming simplicity, one can appeal to the 
classification of Lie algebras. In this way, using the results of 
V.~Kac on automorphisms of Lie algebras (cf.~\cite{Helgason}),
an explicit list of all the generalised symmetric spaces of compact 
simple simply connected Lie groups is given in~\cite{T0}. 

Even when $G$ is not simply connected, 
by~\cite{T0} one has a list of possible
generalised symmetric spaces given by the classification of the 
generalised symmetric Lie algebras, or, equivalently, the 
generalised symmetric spaces of the simply 
connected groups.

From the classification one concludes that generalised symmetric 
spaces $G/H$ with $G$ compact, simple and simply connected and
with $\rk H <\rk G$ occur only for the groups $SU(n)$, $Spin(2n)$ 
and $E_{6}$, compare~\cite{T}.

For generalised symmetric spaces, the application of Cartan's theorem 
is made possible by describing the inclusion of the maximal abelian 
subalgebra of the subgroup $H$ into the maximal abelian subalgebra of $G$.
More precisely, in~\cite{T} the second author gave an explicit formula
expressing, for an arbitrary automorphism $\Theta$, a basis of 
$\TT^{\Theta}$ through a basis of $\TT$.
This formula makes it possible to proceed to explicit calculations of
the map $\rho^{*}(H,G)$ for the generalised symmetric spaces.
The first result is:
\begin{thm}[\cite{T}]\label{gss1}
All generalised symmetric spaces of simple compact Lie groups
are of Cartan type.
\end{thm}
Because of Theorem~\ref{first}, calculations of the cohomology are 
of interest only in the cases where $\rk H < \rk G$. By the classification, 
in almost all cases there is then only one possibility for $\rk H$, 
the exception being $G=Spin(8)$. For these spaces one has:
\begin{thm}[\cite{T}]\label{gss2}
The real cohomology algebra of a generalised symmetric space $G/H$ with
$\rk H < \rk G$ is as follows:
\begin{enumerate}
\item If $G=SU(l+1)$, then
\begin{itemize}    
\item for $l=2n$, $n\geq 1$
\[ 
H^{*}(G/H)\cong 
({\R}[\s]^{W_H}/\langle\rho^{*}(H,G)\sigma_{j}(x_1^2,\ldots,x_n^2)\rangle)
\otimes \wedge (z_3,\ldots,z_{2n+1})  \ ,
\]
\item for $l=2n-1$, $n\geq 3$
\[ 
H^{*}(G/H)\cong 
({\R}[\s]^{W_H}/\langle\rho^{*}(H,G)\sigma_{j}(x_1^2,\ldots,x_n^2)\rangle)
\otimes \wedge (z_3,\ldots,z_{2n-1}) \ ,
\]
\end{itemize}
where $\sigma_{i}$ are the elementary symmetric functions, the $z_{i}$ are
universal transgressive generators of $H^{*}(G)$ corresponding
to $\sigma _{i}$ by transgression in the universal bundle for $G$.
\item If $G=Spin(2n+2)$, $n\geq 2$, and $\rk H = n$, then
\begin{equation*}
H^{*}(G/H)\cong 
({\R}[\s]^{W_H}/\langle \rho^{*}(H,G)\sigma_j(x_1^2,\ldots ,x_n^2) \rangle ) 
\otimes \wedge (z_{n+1}) \ ,
\end{equation*}
where $z_{n+1}$ is a universal transgressive generator of
$H^{*}(G)$ corresponding to $x_1\ldots x_n$.
\item If $G=Spin(8)$ and $\rk H =2$, then
\begin{alignat*}{1}
&H^{*}(G/H)\cong \\
&({\R}[\s]^{W_H}/\langle \rho^{*}(H,G)\sigma_1(x_1^2,x_2^2,x_3^2)
\rho^{*}(H,G)\sigma_3(x_1^2,x_2^2,x_3^2)\rangle )\otimes \wedge (z_2,z_4) \ ,
\end{alignat*}
where $z_2$, $z_4$ are universal transgressive generators of $H^{*}(G)$
corresponding to $\sigma_{2}$ and $x_1x_2x_3x_4$ respectively.
\item If $G=E_{6}$, then
\begin{alignat*}{1}
&H^{*}(G/H)\cong \\
&({\R}[\s]^{W_H}/\langle \rho^{*}(H,E_6)I_2,
\rho^{*}(H,E_6)I_6, \rho^{*}(H,E_6)I_8, \rho^{*}(H,E_6)I_{12}\rangle
)\otimes\wedge (z_5,z_9) \ ,
\end{alignat*}
where $I_2$, $I_6$, $I_8$, $I_{12}$ are the generators of the 
Weyl invariants given in~\cite{Takeuchi}, and $z_5$, $z_9$ are 
universal transgressive generators of 
$H^{*}(E_{6})$ corresponding to  $I_5$, $I_9$.
\end{enumerate}
\end{thm}

The above theorem, together with Lemma~\ref{l:colapses}, implies the 
following. 
\begin{lem}\label{l:surj}
For any two generalised symmetric spaces $G/H$ and $G/L$ of a simple 
compact Lie group $G$ such that $\rk H = \rk L$ and $H\subset L$, the fibration 
$G/H \rightarrow G/L$ with fiber $L/H$ has the property that restriction 
to the fiber is a surjection in cohomology. 
\end{lem} 

\section{Formality in the sense of Sullivan}\label{s:formal}

We will show in this section that for generalised symmetric spaces 
formality in the sense of Sullivan is an immediate consequence of 
their cohomology structure.

Recall that a differentiable manifold is said to be formal in the sense 
of Sullivan if its de Rham algebra of differential forms and its 
cohomology algebra endowed with the zero differential are weakly equivalent, 
meaning that they can be connected by a sequence of 
quasi-isomorphisms, compare~\cite{tokyo} or~\cite{GM}.

In the mid-1970s it became clear~\cite{lehm} that the 
Cartan algebra of a homogeneous space contains more information on its 
topology than that given by Cartan's theorem. More precisely, it turned 
out that the algebra of differential forms on a homogeneous space is weakly 
equivalent to its Cartan algebra. Thus, for homogeneous spaces formality is 
equivalent to formality of its Cartan algebra~\cite{Onishchik}. 
However, formality of the Cartan algebra can be
described in terms of its cohomology:
\begin{thm}\label{t:form}
For the Cartan algebra $(C,d)$ of a compact homogeneous space $G/H$ with
$\rk G = n$ and $\rk H = k$ the following conditions are equivalent:
\begin{enumerate}
\item $H^{*}(C) = H^{*}(BH)/\langle \rho^{*}(H, G)H^{*}(BG)\rangle \otimes
\wedge (z_{k+1},\ldots,z_{n})$,
\item $(C,d)$ is formal.
\end{enumerate}
\end{thm}
A more general statement on the formality of Cartan algebras can be
found in~\cite{GHV,Onishchik} in the context of formal algebras.

\begin{rem}
Theorem~\ref{t:form} implies that a compact homogeneous space is formal if
and only if it is of Cartan type. It follows immediately that all
homogeneous 
spaces $G/H$ with $\rk G = \rk H$ are formal.
\end{rem}

\begin{rem}
From the cohomology calculations in~\cite{Borel} and~\cite{Takeuchi} 
described in subsection~\ref{ss:symm} above it follows that all
symmetric spaces are of Cartan type. Together with Theorem~\ref{t:form} 
this shows that compact symmetric spaces are formal. This is usually
proved by showing that they are geometrically formal, as 
in~\cite{DFN}, but the cohomological proof seems to be closer in 
spirit to Sullivan's theory of minimal models.
\end{rem}

Combining Theorem~\ref{t:form} with Theorem~\ref{gss1} we obtain:
\begin{thm}\label{t:Sform}
All generalised symmetric spaces of simple compact Lie groups are formal
in the sense of Sullivan.
\end{thm}

\begin{rem}
Theorem~\ref{t:Sform} extends partial results due to Duma\'nska-Malyszko, 
Stepie\'n and Tralle~\cite{DST}.
We were recently informed by A.~Tralle that a different proof
of Theorem~\ref{t:Sform} has been given by Z.~Stepie\'n~\cite{Stepien}. \end{rem}

\section{The cohomology structure of flag manifolds}\label{s:flag}

Our proof that the non-symmetric flag manifolds are not geometrically 
formal requires detailed, explicit, information about the generators 
and relations of their cohomology rings. This section provides these 
details, starting from the theorems of Section~\ref{s:gsscoh}. We 
originally obtained the formulas for the relations using computer 
calculations with Groebner basis algorithms. Having found the 
formulas, they are then easy to prove by elementary arguments not 
invoking Groebner bases.

The first case to consider is that of the classical flag manifolds 
$SU(n+1)/T^n$. From Theorem~\ref{first} we have:
\[ 
H^{*}(SU(n+1)/T^n)\cong {\R}[x_0,...,x_n]/\langle
S^{+}(x_0,...,x_n)\rangle \ , 
\]
where the $S^{+}(x_0,...,x_n)$ are the symmetric functions of positive
degrees.
Here is an explicit form:
\begin{prop}\label{p:basis}
The classes represented by
\begin{equation}\label{gener}
x_{1}^{\alpha_1}x_{2}^{\alpha_2}\ldots x_{n}^{\alpha _n} \ ,
\qquad 0\leq \alpha_{i}\leq i \ , \qquad 1\leq i\leq n
\end{equation}
form a basis for the cohomology of $SU(n+1)/T^n$ as a vector space.
Multiplicative relations between the $x_1,\ldots,x_n$ are given by:
\begin{equation}\label{rel}
\sum _{i_{1}+\ldots+i_{p} =
n-p+2}x_{n-p+1}^{i_p}x_{n-p+2}^{i_{p-1}}\ldots x_{n-1}^{i_2}x_{n}^{i_1}
= 0, \qquad 1\leq p\leq n \ .
\end{equation}
\end{prop}
\begin{proof}
Define $s_{k}(x_{0},\ldots,x_{n})$ to be the sum of all the 
monomials in the $x_{i}$ which are homogeneous of degree $k$. This is 
clearly a symmetric polynomial.    
The relations~\eqref{rel} that we have to prove amount to
\begin{equation}
s_{k}(x_{k-1},\ldots,x_{n})=0,\qquad 2\leq k\leq n+1 \ .
\end{equation}
We shall prove more, namely that 
\begin{equation}\label{eq:rell}
s_{m}(x_{k-1},\ldots,x_{n})=0,\qquad 2\leq k\leq n+1 \ ,
\end{equation}
holds for all $m\geq k$.

First we prove that $s_{m}(x_{1},\ldots,x_{n})=0$ for all
$m\geq 2$. We know
$$
0=s_{m}(x_{0},\ldots,x_{n})=x_{0}s_{m-1}(x_{0},\ldots,x_{n})
+s_{m}(x_{1},\ldots,x_{n}) \ ,
$$
because all the symmetric functions in $x_{0},\ldots,x_{n}$ vanish.
Thus the vanishing of $s_{m-1}(x_{0},\ldots,x_{n})$ implies the 
vanishing of $s_{m}(x_{1},\ldots,x_{n})$. 

We now prove~\eqref{eq:rell} by induction on $k$. The case $k=2$ is 
what we just proved. Suppose we have proved the statement up to $k$. 
Then consider
$$
s_{m+1}(x_{k-1},\ldots,x_{n})=x_{k-1}s_{m}(x_{k-1},\ldots,x_{n})
+s_{m+1}(x_{k},\ldots,x_{n}) \ .
$$
As soon as $m\geq k$ both the left hand side and the first summand 
on the right hand side vanish by the induction hypothesis. Therefore the 
second summand on the right also vanishes, which is what is to be 
proved in the inductive step.

Having proved the multiplicative relations, it remains to prove the 
statement about the vector space basis of the cohomology. This can be 
proved by induction on the degree. A vector space basis for $H^{2}$ 
is given by $x_{1},\ldots,x_{n}$. Suppose we have proved the 
statement up to degree $2k$. Now $H^{2k+2}$ is linearly generated by 
all homogeneous monomials of degree $k+1$ in the $x_{1},\dots,x_{n}$.
However, by induction there are linear relations expressing 
$x_{1}^{2}$ as a linear combination of monomials containing $x_{1}$ 
at most linearly, expressing $x_{2}^{3}$ as a linear 
combination of monomials containing $x_{2}$ at most in the second power,
and so on up to $x_{k-1}^{k}$. We also have a new relation in this 
degree, namely $s_{k+1}(x_{k},\ldots,x_{n})=0$. This allows us 
to replace $x_{k}^{k+1}$ by a linear combination of monomials 
containing only smaller powers of $x_{k}$.

We now have eliminated all monomials not listed in~\eqref{gener}.
The remaining ones must be linearly independent because their number 
in each degree is seen to equal the respective Betti number by 
inspection of the Poincar\'e polynomial.
\end{proof}

Next we consider the flag manifolds $Spin(2n+1)/T^n = SO(2n+1)/T^n$ 
and $Sp(n)/T^n$. Theorem~\ref{first} implies
\[
H^{*}(Spin(2n+1)/T^n) = H^{*}(Sp(n)/T^n)\cong {\R}[x_1,...,x_n]/\langle
S^{+}(x_1^2,...,x_n^2)\rangle \ .
\]
We can use the same argument as in the proof of the previous 
proposition to obtain: 
\begin{prop}\label{p:basis1}
The classes represented by
\begin{equation}
x_{1}^{\alpha_1}x_{2}^{\alpha_2}\ldots x_{n}^{\alpha _n} \ ,
\qquad 0\leq \alpha_{i}\leq 2i-1 \ , \qquad 1\leq i\leq n
\label{gener1}
\end{equation}
form a vector space basis for the cohomology of $Spin(2n+1)/T^n$ and 
of $Sp(n)/T^n$. Multiplicative relations between the 
$x_1,\ldots,x_n$ are given by:
\begin{equation}\label{rel1}
\sum _{i_{1}+\ldots+i_{p} =
n-p+1}x_{n-p+1}^{2i_p}x_{n-p+2}^{2i_{p-1}}\ldots x_{n-1}^{2i_2}x_{n}^{2i_1}
= 0, \qquad 1\leq p\leq n \ .
\end{equation}
\end{prop}

Finally, we consider $Spin(2n)/T^n = SO(2n)/T^n$. 
In this case Theorem~\ref{first} gives:
\[
H^{*}(Spin(2n)/T^n)\cong {\R}[x_1,\ldots ,x_n]/
\langle S^{+}(x_1^2,\ldots ,x_n^2), \ x_1\ldots x_n \rangle \ .
\]
More explicitly:
\begin{prop}\label{p:basis2}
A vector space basis for the cohomology of $Spin(2n)/T^n$ is given by
\begin{equation}
x_{1}^{\alpha_1}x_{2}^{\alpha_2}\ldots x_{n}^{\alpha_n} \ ,
\label{gener2}
\end{equation}
with the coefficients $\alpha _i$ satisfying: $0\leq \alpha _{i}\leq 2i-1$ 
for $1\leq i\leq n-1$, $ 0\leq \alpha_{n}\leq 2n-2$ and $\alpha _{i}=2i-1$ 
implies $\alpha _{i+1}\ldots \alpha_{n}=0$.
 
Multiplicative relations between the $x_1,\ldots ,x_n$ are given by the 
formulas
\begin{equation}\label{rel2}
\sum_{i_1+\ldots +i_p = n-p+1} x_{n-p+1}^{2i_p}x_{n-p+2}^{2i_{p-1}}\ldots 
x_{n-1}^{2i_2}x_{n}^{2i_1} = 0 \ , \quad 1\leq p\leq n \ ,
\end{equation}
\begin{equation}\label{rel3}
\sum_{i_1+\ldots +i_p = n} x_{n-p+1}^{2i_p-1}x_{n-p+2}^{2i_{p-1}-1}\ldots 
x_{n-1}^{2i_2-1}x_{n}^{2i_1-1} = 0 \ , \quad 1\leq p\leq n \ .
\end{equation}
\end{prop}
\begin{proof}
To prove the relations~\eqref{rel2} we can proceed as in the proof of 
Propositions~\ref{p:basis} and~\ref{p:basis1}. 

To prove the relations~\eqref{rel3} we will proceed by backward induction
on $p$.
For $p=n$ the left hand side is $x_1\ldots x_n$, which obviously vanishes. 
Now suppose we have proved the statement down to $p-1\leq n$. Consider the 
relation
\[
 \sum_{i_1+\ldots +i_{p+1} = n-p} x_{n-p}^{2i_{p+1}}x_{n-p+1}^{2i_{p}}
 \ldots x_{n-1}^{2i_2}x_{n}^{2i_1} = 0 
\]
which was proved already. Multiplying it by $x_{n-p+1}\ldots x_{n}$ and 
splitting the resulting sum into two sums corresponding to the cases 
$i_{p+1}\neq 0$ and $i_{p+1}=0$ we get
\begin{alignat*}{1}
x_{n-p} \sum_{i_1+\ldots +i_{p+1}=n-p}
&x_{n-p}^{2i_{p+1}-1}x_{n-p+1}^{2i_p+1}\ldots x_n^{2i_1+1} \\
+ &\sum_{i_1+\ldots 
+i_{p}=n-p}x_{n-p+1}^{2i_{p}+1}x_{n-p+2}^{2i_{p-1}+1}
\ldots x_n^{2i_1+1} = 0 \ .
\end{alignat*}
The first sum vanishes by the induction hypothesis.
Therefore the second sum vanishes, which is what we need to prove in 
the inductive step.

To prove the statement about the vector space basis of the cohomology 
we can use the same argument as in the proof of Proposition~\ref{p:basis}.
\end{proof}

\begin{rem}
We originally obtained the formulas discussed above using Groebner basis
algorithms. It can be proved that the
polynomials defining the relations~\eqref{rel}, \eqref{rel1} 
and~\eqref{rel2}, \eqref{rel3} give Groebner bases for the ideals 
$\langle S^{+}(x_0,\ldots ,x_n)\rangle$, 
$\langle S^{+}(x_1^2,\ldots ,x_n^2)\rangle$ and 
$\langle S^{+}(x_1^2,\ldots ,x_n^2), x_1\ldots x_n\rangle$ respectively.
This also implies that the polynomials given by~\eqref{gener}, \eqref{gener1} and~\eqref{gener2} form vector space bases for the 
corresponding cohomology algebras.
\end{rem}
 
\section{Failure of geometric formality}\label{s:geoform}

In this section we prove that various generalised symmetric
spaces $G/H$ are not geometrically formal in the sense of~\cite{K},
i.~e.~that they do not admit Riemannian metrics for which all products 
of harmonic forms are harmonic. Note that we do not assume that the 
metrics are $G$-invariant.

The simplest result, which however illustrates a main part of the 
argument for the flag manifolds as well, is the following theorem.
\begin{thm}\label{t:G2}
The $6$-symmetric space $G_{2}/T^{2}$ is not geometrically formal.
\end{thm}
\begin{proof}
The real cohomology of $X=G_{2}/T^{2}$ was calculated by
Borel~\cite{Borel}; we use the presentation in~\cite{BS}. There are
two linearly independent generators $x$ and $y\in H^{2}(X,\R)$, which
satisfy the relations
\begin{equation}\label{eq:G2rel1}
x^{2}+3xy+3y^{2}=0
\end{equation}
and  
\begin{equation}\label{eq:G2rel2}
x^{6}=y^{6}=0 \ .
\end{equation}
On the other hand, $xy^{5}$ generates the top-dimensional cohomology
$H^{12}(X,\R)$.

Suppose that $X$ admits a formal Riemannian metric. By an obvious abuse
of notation, we denote by $x$ and $y$ the harmonic representatives of
the above cohomology classes. Then the above relations for $x$ and $y$
hold at the level of differential forms. In particular $x\wedge y^{5}$
is a volume form on $X$.

On the other hand, it follows from~\eqref{eq:G2rel2} that both
kernel distributions
\[ 
N_{x} = \{ v\in TM \ \ \vert \ \  i_{v}x= 0 \}
\]
\[ 
N_{y} = \{ w\in TM \ \ \vert \ \  i_{w}y= 0 \}
\]
have rank at least $2$. Therefore, we can locally choose linearly
independent vector fields $v\in N_{x}$ and $w\in N_{y}$. It follows
from~\eqref{eq:G2rel1} that $i_{w}x\wedge i_{v}y=0$. But this
implies $i_{v}i_{w}(x\wedge y^{5})=0$, contradicting the fact that
$x\wedge y^{5}$ is a volume form.
\end{proof}

We now consider the flag manifolds.
\begin{thm}\label{t:nonformal}
For all $n\geq 2$
the classical flag  manifolds  $SU(n+1)/T^n$ are not geometrically formal.
\end{thm}
By the results of~\cite{T0}, $SU(n+1)/T^n$ is a generalised symmetric 
space of order $n+1$. For $n=1$ it is the symmetric space
$S^{2}$, which is of course geometrically formal.

The proof of Theorem~\ref{t:nonformal} uses the same idea as that of 
Theorem~\ref{t:G2}, together with induction over $n$. To carry this 
out, we need the explicit relations from Proposition~\ref{p:basis}.
First, we will prove the following lemma.
\begin{lem}\label{l:vanish}
Let $M$ be a differentiable manifold of dimension $n^2 + n$, with
$n\geq 2$. Suppose there are $n$ closed two-forms $x_1,\ldots,x_n$ on
$M$ satisfying relations~\eqref{rel}. Then $x_1\wedge
x_2^2\wedge\ldots\wedge x_n^n$ vanishes identically. In particular,
it is not a volume form on $M$.
\end{lem}
\begin{proof}
Note that the first relation in~\eqref{rel}, for $p = 1$, gives 
$x_{n}^{n+1} = 0$. Using this, the second relation, for $p=2$, 
implies $x_{n-1}^{n+1} = 0$.

Assume that, under the assumptions of the lemma,
$x_1\wedge x_2^2\wedge\ldots\wedge x_n^n$ is a volume form on some 
open subset $U\subset M$. Then
$x_n^{n+1}=0$, but $x_{n}^{n}$ is not zero on $U$. Thus $x_{n}$ is of
constant rank $2n$ in $U$.
By the same argument, $x_{n-1}$ is of constant rank equal to $2n-2$ or
$2n$ in $U$. The kernel distributions
\[ 
N_{x_n} = \{ w\in TM \ \ \vert \ \  i_{w}x_{n}= 0 \}
\]
\[ 
N_{x_{n-1}} = \{ v\in TM \ \ \vert \ \  i_{v}x_{n-1}= 0 \}
\]
are of ranks $\rk (N_{x_n})=n^{2}-n$ and $\rk (N_{x_{n-1}})=n^2 - n + 2$
or $n^2 - n$. As $x_{n}$ and $x_{n-1}$ are closed, the
Frobenius theorem implies that the kernel distributions are integrable.

We proceed by induction on $n$.
For $n =2$ formula~\eqref{rel}
gives the following  relation between $x_1$ and $x_2$:
\begin{equation}\label{f3}
    x_1^2 + x_{1}\wedge x_{2} + x_2^2 = 0 \ .
    \end{equation}
From the above discussion, the ranks of $N_{x_1}$ and $N_{x_2}$ are
$\geq 2$. Thus locally there are linearly independent vectors
$v\in N_{x_1}$ and $w\in N_{x_2}$. From~\eqref{f3} it follows that
$i_{w}x_{1}\wedge i_{v}x_{2} = 0$,  which implies
$i_{w}i_{v}(x_{1}\wedge x_{2}^2) = 0$. This implies that
$x_{1}\wedge x_{2}^2$ can not be a volume form anywhere, and therefore
vanishes identically.

Assume that the statement holds for $n - 1\geq 2$. Let us consider a
manifold $M$ of dimension $n^2 + n$ and let $x_1,\ldots,x_n$ be forms
on $M$ satisfying~\eqref{rel}, such that
$x_{1}\wedge x_{2}^2\wedge\ldots\wedge x_{n}^n$ is a volume form on 
some open subset $U\subset M$.
Since $N_{x_{n}}$ is integrable, it defines a foliation. Let $N$ be a
leaf of this foliation. Then $N$ is of dimension
$n^2 - n = (n - 1)^2 + (n - 1)$, and the forms $x_1,\ldots,x_{n-1}$
restricted to $N$ satisfy relations~\eqref{rel} and
$x_{1}\wedge x_{2}^2\wedge\ldots\wedge x_{n-1}^{n-1}$ is a volume
form on $N$. This contradicts the induction hypothesis.
\end{proof}

\begin{proof}[Proof of Theorem~\ref{t:nonformal}]
Assume that the flag manifold $SU(n+1)/T^n$  is geometrically formal,
that is, there is a metric for which all products of harmonic forms
are harmonic. If for the classes $x_1,\ldots,x_n$ we choose their
harmonic representatives (denoted by the same letters), geometric
formality implies that the relations~\eqref{rel} hold at the level of
differential forms. 
The dimension of
$SU(n+1)/T^n$ is $n^2 + n$, and from formula~\eqref{gener} we see that
$x_{1}\wedge x_{2}^2\wedge\ldots\wedge x_{n}^{n}$ is a volume form on
$SU(n+1)/T^n$. This gives a contradiction with the above lemma.
\end{proof}

\begin{thm}\label{t:nonformal1}
For all $n\geq 2$ the flag manifolds $Spin(2n+1)/T^n$ and $Sp(n)/T^n$ 
are not geometrically formal.
\end{thm}
By~\cite{T0}, these spaces are generalised symmetric of order $2n$. 
For $n=1$ we again obtain the $2$-sphere.

As before, we first prove a lemma about closed forms satisfying certain 
relations.
\begin{lem}\label{l:vanish1}
Let $M$ be a smooth manifold of dimension $2n^2$, with
$n\geq 2$. Suppose there are $n$ closed two-forms $x_1,\ldots,x_n$ on
$M$ satisfying the relations~\eqref{rel1}. Then $x_1\wedge
x_2^3\wedge\ldots\wedge x_n^{2n-1}$ vanishes identically. In 
particular, it is not a volume form on $M$.
\end{lem}
\begin{proof}
We proceed as in the proof of Lemma~\ref{l:vanish}.
The first relation in~\eqref{rel1}, for $p = 1$, gives 
$x_{n}^{2n} = 0$. Using this, the second relation, for $p=2$, 
implies $x_{n-1}^{2n} = 0$.

If we assume that $x_1\wedge x_2^3\wedge\ldots\wedge x_n^{2n-1}$ is a 
volume form on some open subset $U\subset M$, then in $U$ we conclude 
that $x_{n}$ is of
constant rank $2(2n-1)$ and  $x_{n-1}$ is of constant rank equal to $2(2n-3)$, 
$2(2n-2)$ or $2(2n-1)$. So, the kernel distributions
\[ 
N_{x_n} = \{ w\in TM \ \ \vert \ \  i_{w}x_{n}= 0 \}
\]
\[ 
N_{x_{n-1}} = \{ v\in TM \ \ \vert \ \  i_{v}x_{n-1}= 0 \}
\]
are of ranks $\rk (N_{x_n})=2n^{2} - 4n + 2$ and 
$\rk (N_{x_{n-1}})=2n^2 - 4n + 6$, $2n^2-4n+4$ or $2n^2 - 4n + 2$, and are 
integrable.

We proceed by induction on $n$.
For $n =2$ formula~\eqref{rel1}
gives the following  relation between $x_1$ and $x_2$:
\begin{equation}\label{f4}
    x_1^2 +  x_2^2 = 0 \ .
    \end{equation}
From the above discussion , the rank of  $N_{x_1}$ is $\geq 2$, thus  
locally there is a non-zero vector field $v\in N_{x_1}$.  From~\eqref{f4} 
it follows that $x_{2}\wedge i_{v}x_{2} = 0$,  which implies
$i_{v}(x_{1}\wedge x_{2}^3) = 0$. This implies that
$x_{1}\wedge x_{2}^3$ can not be a volume form.

Assume that the statement holds for $n - 1\geq 2$. Let us consider a
manifold $M$ of dimension $2n^2$ and let $x_1,\ldots,x_n$ be forms
on $M$ satisfying~\eqref{rel1}, such that
$x_{1}\wedge x_{2}^3\wedge\ldots\wedge x_{n}^{2n-1}$ is a volume form on $M$.
Since $N_{x_{n}}$ is integrable, it defines a foliation. Let $N$ be a
leaf of this foliation. Then $N$ is of dimension
$2n^{2} - 4n + 2 = 2(n - 1)^2 $, and the forms $x_1,\ldots,x_{n-1}$
restricted to $N$ satisfy relations~\eqref{rel1} and
$x_{1}\wedge x_{2}^3\wedge\ldots\wedge x_{n-1}^{2n-3}$ is a volume
form on $N$. This contradicts the induction hypothesis.
\end{proof}

\begin{proof}[Proof of Theorem~\ref{t:nonformal1}]
Assume that $M=Spin(2n+1)/T^{n}$ or $Sp(n)/T^n$ is geometrically formal. 
If for the classes $x_1,\ldots,x_n$ we choose their harmonic 
representatives, geometric formality 
implies that the relations~\eqref{rel1} hold at the level of
differential forms. 
The dimension of $M$ is $2n^2$, and from formula~\eqref{gener1} 
we see that $x_{1}\wedge x_{2}^{3}\wedge\ldots\wedge x_{n}^{2n-1}$ 
is a volume form on $M$. This contradicts the above lemma.
\end{proof}

\begin{thm}\label{t:nonformal2}
For all $n\geq 4$ the flag manifolds $Spin(2n)/T^n$ are not 
geometrically formal.
\end{thm}
By~\cite{T0}, $Spin(2n)/T^n$ is generalised symmetric of order $2n-2$.
For $n=2$ we obtain the symmetric space $S^{2}\times S^{2}$. 
For $n=3$ we obtain $SU(4)/T^3$, which by Theorem~\ref{t:nonformal}
is not geometrically formal. 
 
The proof of Theorem~\ref{t:nonformal2} is more complicated than
that of the previous ones, because the cohomology algebra is more complicated.
We shall first prove the following:
\begin{prop}\label{l:vanish2}
Let $M$ be a smooth manifold of dimension $2n^2-2n$, with
$n\geq 4$. Suppose there are $n$ closed two-forms $x_1,\ldots,x_n$ on
$M$ satisfying relations~\eqref{rel2} and ~\eqref{rel3}. Then $x_2^2\wedge
x_3^4\wedge \ldots \wedge x_n^{2n-2}$ vanishes identically. In 
particular, it is not a volume form on $M$.
\end{prop}
Note that for $p=1$ the relation~\eqref{rel3} becomes $x_n^{2n-1}=0$. 
Also from~\eqref{rel3} the following relation can easily be obtained by 
backward induction on $k$:
\begin{equation}\label{addrel}
x_{k}^{2k-1}x_{k+1}^{2k-1}x_{k+2}^{2k+1}\ldots 
x_{n-1}^{2n-5}x_{n}^{2n-3} = 0 \ , \quad 2\leq k\leq n-1 \ .
\end{equation}

Put $M_{n+1} = M$ and recursively define $M_{k}$ to be a leaf of the kernel 
foliation of $x_{k}$ restricted to $M_{k+1}$, for all $2\leq k\leq n$. 
\begin{lem}
Let $M$ be a smooth manifold of dimension $2n^2-2n$, with $n\geq 4$. 
Suppose there are $n$ closed $2$-forms satisfying relations~\eqref{rel3} such that
$x_2^2\wedge \ldots \wedge x_{n}^{2n-2}$ is a volume form on $M$. Then 
$x_{k}^{2k-1}$ vanishes identically on $M_{k+1}$ and 
$x_2^2\wedge \ldots \wedge x_{k-1}^{2(k-1)-2}$ is a volume form on $M_{k}$.
\end{lem}
\begin{proof}
We will proceed by backward induction on $k$.
To prove this for $k=n$ note that $x_{n}^{2n-1}=0$ on $M$ and the 
assumption that $x_2^2\wedge \ldots \wedge x_{k-1}^{2n-2}$ is a 
volume form on $M$ implies that $x_{n}$ has constant rank equal 
to $4n-4$, so, being a leaf of its kernel foliation, 
$M_{n}$ has dimension $2(n-1)^2-2(n-1)$ and 
$x_2^2\wedge \ldots \wedge x_{k-1}^{2(n-1)-2}$ 
is a volume form on $M_n$.

Assume that the lemma has been proved for all $k+1\geq 4$. Since 
$x_2^2\wedge\ldots\wedge x_{k}^{2k-2}$ is a volume form on $M_{k+1}$ 
we conclude that $\dim M_{k+2}-\dim M_{k+1}=4k$, for all $k+1\geq 4$.
As $M_k\subset M_{k+1}$, denote by $D_k$ a distribution complementary 
to $TM_k$ in $TM_{k+1}$. Relation~\eqref{addrel} implies that the form
$x_{k}^{2k-1}\wedge x_{k+1}^{2k-1}\wedge\ldots\wedge x_{n}^{2n-3}$ 
vanishes identically on $M$. If we evaluate this form on $2(2k-1)$ 
vectors from $TM_{k+1}$, $2(2k-1)$ vectors from $D_{k+1}$ on which 
$x_{k+1}^{2k-1}$ does not vanish, $2(2k+1)$ vectors form 
$D_{k+2}$ on which $x_{k+2}^{2k+1}$ does not vanish, and so on, and
finally $2(2n-3)$ vectors from $D_{n}$ on which $x_n^{2n-3}$ does
not vanish, 
we conclude that $x_{k}^{2k-1}$ vanishes identically on $M_{k+1}$. Since 
$\dim D_{k+1}=4k$, for $k+1\geq 4$, the choice of such a vectors 
is always possible.

Since $x_2^2\wedge\ldots\wedge x_{k}^{2k-2}$ is a volume form on 
$M_{k+1}$ it follows that $x_{k}$ restricted to $M_{k+1}$ has constant 
rank equal to $2k-2$. Thus, $\dim M_{k}=2(k-1)^2-2(k-1)$, so 
$x_2^2\wedge\ldots\wedge x_{k-1}^{2(k-1)-2}$ is a volume form on $M_k$.
\end{proof} 

\begin{proof}[Proof of Proposition~\ref{l:vanish2}] 
Assume that under the conditions given in the proposition, 
$x_2^2\wedge x_3^4\ldots\wedge x_n^{2n-2}$ is a volume form on $M$. 
Then the above lemma implies that we have the following situation: 
a manifold $M$ of dimension $2n^2-2n$ and $n$ closed $2$-forms 
satisfying relations~\eqref{rel2} such that $x_k$ restricted to $M_{k+1}$ 
has constant rank equal to $2k-2$ and 
$x_2^2\wedge\ldots\wedge x_{k-1}^{2(k-1)-2}$ is a volume form on $M_k$. 
Note that the forms $x_2,\ldots,x_{k}$ satisfy the relations~\eqref{rel2} 
on $M_{k+1}$. 

We prove by induction on $n$ that this situation gives a contradiction.

For $n=4$ we have that $x_2^2\wedge x_3^4$ is a volume on $M_4$ and
\eqref{rel3} implies that on $M_4$ we have following relations:
\[
x_2^4+x_2^2\wedge x_3^2+x_3^4=0,\;\; x_3^6=0 \ .
\]
As in the proofs of Lemmas~\ref{l:vanish} and~\ref{l:vanish1} this 
gives a contradiction.

Let us assume that the statement holds for all $n-1\geq 4$. 
Consider a manifold $M$ of dimension 
$2n^2-2n$ and assume we have $2$-forms $x_1,\ldots ,x_n$ satisfying above 
conditions. Then, obviously, we have on $M_n$ the same situation with 
$n-1$ two-forms 
$x_1,\ldots,x_{n-1}$ giving the contradiction.
\end{proof}

\begin{proof}[Proof of Theorem~\ref{t:nonformal2}]
Assume that $M=Spin(2n)/T^{n}$, $n\geq 4$  is geometrically formal. 
If for the classes $x_1,\ldots,x_n$ we choose their harmonic 
representatives (denoted by the same letters), geometric formality 
implies that the relations~\eqref{rel2} and~\eqref{rel3} hold at the level of
differential forms. 
The dimension of $M$ is $2n^2-2n$, and from formula~\eqref{gener2} 
we see that $x_{2}^2\wedge x_{3}^{4}\wedge\ldots\wedge x_{n}^{2n-2}$ 
is a volume form on $M$. This contradicts Proposition~\ref{l:vanish2}.
\end{proof}
  
So far we have only considered homogeneous spaces $G/H$ where $G$ and 
$H$ have equal rank. There is a generalisation of these arguments to 
some generalised symmetric spaces $G/H$ with $\rk G > \rk H$. 
The simplest case is the following:
\begin{thm}\label{t:Spin8}
The $12$-symmetric space $X=Spin(8)/T^{2}$ is not geometrically
formal. 
\end{thm}
\begin{proof}
By Theorem~\ref{gss2} the cohomology algebra of $X$ is the tensor product
of a polynomial algebra $P$, which is the cohomology algebra of
$G_{2}/T^{2}$, and an exterior algebra $E$, which is the cohomology
algebra of $S^{7}\times S^{7}$.

The inclusions $T^{2}\subset G_{2}\subset Spin(8)$ induce a fibration
$X\longrightarrow Z=Spin(8)/G_{2}$ with fiber
$Y=G_{2}/T^{2}$. As the base and the total space are generalised
symmetric spaces, Lemma~\ref{l:surj} implies that all cohomology classes on $Y$ 
are restrictions of classes on $X$.

We shall use the basis for $P$ used in the proof of Theorem~\ref{t:G2}.
Assume that $X$ is geometrically formal and identify all the elements
of the cohomology algebra with their harmonic representatives. Then
the harmonic $2$-forms $x$ and $y$ on $X$ satisfy $x^{6}=0=y^{6}$,
but $x^{5}\neq 0\neq y^{5}$, and therefore have kernels of rank
$\dim (X)-10 = 16$. As the codimension of the fiber $Y$ in $X$ is
$14$, it follows that the restrictions of $x$ and $y$ to $Y$ have
kernels of rank at least $2$ everywhere.

Thus at every point of a fiber we can find linearly independent
local vector fields $v$ and $w$ contained in the kernels of $x$ and
$y$ respectively. As the restrictions of $x$ and $y$ to $Y$ satisfy
relation~\eqref{eq:G2rel1}, we conclude $i_{v}i_{w}(x\wedge y^{5})=0$
as in the proof of Theorem~\ref{t:G2}. This shows that the restriction
of $x\wedge y^{5}$ to $Y$ vanishes identically. This contradicts the
fact that restriction to $Y$ is surjective in cohomology.
\end{proof}

Using the theorems about the flag manifolds in a similar way, we also 
obtain:
\begin{thm}\label{t:nonf}
The following generalised symmetric spaces are not geometrically formal:
\begin{enumerate}
\item $SU(2n+1)/T^n$, for $n\geq 2$,
\item $SU(2n)/T^n$, for $n\geq 3$,
\item $Spin(2n+2)/T^n$, for $n\geq 2$.
\end{enumerate}
\end{thm}
By~\cite{T0}, these are indeed generalised symmetric spaces of order 
$4n+2$, $4n-2$ and $2n+2$ respectively. We could consider $n=2$ in
the second case, this would give $SU(4)/T^2$ which is the same as 
$Spin(6)/T^2$ contained in the third case. 
\begin{proof}
To prove the first statement, let us consider the fibration 
$SU(2n+1)/T^n\lra SU(2n+1)/SO(2n+1)$ with fiber $SO(2n+1)/T^n = 
Spin(2n+1)/T^n$. The base is a symmetric space, so Lemma~\ref{l:surj} 
shows that the restriction to the fiber is 
surjective in cohomology. Theorem~\ref{gss2} implies that
\[
H^{*}(SU(2n+1)/T^n)\cong {\R}[x_1,\ldots ,x_n]/\langle S^{+}(x_1^2,\ldots ,x_n^2)
\rangle \otimes \wedge (z_3,\ldots ,z_{2n+1}) \ .
\]
Assume that $SU(2n+1)/T^n$ is geometrically formal. For the cohomology 
classes $x_1,\ldots,x_n$ we take their harmonic representatives with 
respect to a formal metric. Then the relations~\eqref{rel1} hold for the
harmonic forms, as forms. If we restrict these forms to the fiber 
$Spin(2n+1)/T^n$, Lemma~\ref{l:vanish1} implies that the form 
$x_1\wedge \ldots \wedge x_{n}^{2n-1}$ vanishes. This contradicts the 
fact that the restriction is surjective in cohomology.

For the second case, we consider the fibration $SU(2n)/T^n\lra SU(2n)/Sp(n)$ 
with fiber $Sp(n)/T^n$, where, as above, restriction to the fiber is
surjective in cohomology. Again Theorem~\ref{gss2} implies
\[
H^{*}(SU(2n)/T^n)\cong {\R}[x_1,\ldots ,x_n]/\langle S^{+}(x_1^2,\ldots 
,x_n^2)\rangle \otimes \wedge (z_3,\ldots ,z_{2n-1}) \
\]
and, as in the first case, 
the assumption of geometric formality for $SU(2n)/T^n$ contradicts 
the fact that the restriction 
to the fiber is surjective in cohomology.

For the third case, we have the fibration $Spin(2n+2)/T^n\lra 
Spin(2n+2)/Spin(2n+1)$ 
with fiber $Spin(2n+1)/T^n$, and we can proceed as above.
\end{proof}
\begin{rem}
Note that if $X$ is the total space of a fibration with fiber $Y$,
there is no reason for the restrictions of harmonic forms on
$X$ to be harmonic on $Y$.
\end{rem}

\begin{rem}
It follows from the classification of generalised symmetric spaces 
$G/H$ in~\cite{T0} that the only such spaces where $H$ is a torus of 
rank $\geq 2$ and $G$ is either $G_{2}$ or a simply connected classical 
simple group are the ones we considered in Theorems~\ref{t:G2}, \ref{t:nonformal}, 
\ref{t:nonformal1}, \ref{t:nonformal2}, \ref{t:Spin8} and~\ref{t:nonf}.
The generalised symmetric spaces of the form $SO(n)/T^{k}$ with $k\geq 
2$ can be treated similarly, using the results of~\cite{T0,T}.
\end{rem} 

\bibliographystyle{amsplain}

\bigskip

\end{document}